\documentclass{amsart}

\newtheorem{theorem}{Theorem}[section]
\newtheorem{lemma}[theorem]{Lemma}
\newtheorem{proposition}[theorem]{Proposition}
\newtheorem{corollary}[theorem]{Corollary}

\theoremstyle{definition}

\newtheorem{example}[theorem]{Example}

\theoremstyle{remark}
\newtheorem{remark}[theorem]{Remark}

\numberwithin{equation}{section}

\begin{document}

\title{Mayberry-Murasugi's formula for links in homology 3-spheres}
\author{Joan Porti}
\address{Departament de Matem{\`a}tiques, Universitat Aut{\`o}noma de Barcelona, 08193 Be\-lla\-ter\-ra, Spain}
\email{porti@mat.uab.es}
\thanks{Partially supported by DGICYT through grant BFM2000-0007}

%    General info

\date{\today}


\begin{abstract}
We prove Mayberry-Murasugi's formula for links in homology 3-spheres,
which was proved before only for links in the 3-sphere. Our proof uses
Franz-Reidemeister torsions.
\end{abstract}

\maketitle

\section{Introduction}

Fox's formula computes the order of the first homology group of a finite
cyclic covering of a knot in $S^3$ from its Alexander polynomial
\cite{Fox}. This formula has been generalized by Mayberry and Murasugi for
finite abelian coverings of links in $S^3$ \cite{MM}.
 Here we give a new proof
of this formula using Franz-Reidemeister torsions, that applies to links
not only in $S^3$ but in homology 3-spheres. Results in this direction
have been obtained in \cite{Sakuma} and \cite{HillmanSakuma}.

Let $M^3$ be a closed three-dimensional homology sphere and
$L\subset M^3$ a smooth link with $\mu$ components
$l_1,\ldots,l_{\mu}$. Its exterior is denoted by $E(L)=M^3-N(L)$.
A finite abelian covering $\hat M^3_{\pi}\to M^3$ branched along
$L$ is given by the kernel of an epimorphism
\[
\pi\!:\pi_1E(L)\to G,
\]
where $G$ is a finite abelian group. The set of representations from $G$
to non-zero complex numbers $\xi\! :G\to \mathbf C^*$ is denoted by
$\hat{G}$, and it is a group isomorphic to $G$, called the Pontrjagin
dual.

We choose meridians $m_1,\dots,m_{\mu}\in H_1(E(L),\mathbb Z)$.
For $\xi\in \hat{G}$, let $L_{\xi}\subseteq L$ be the sublink
consisting of those components $l_i$ such that $\xi(m_i)\neq 1$.
Let $\Delta_{L_\xi}(t_{i_1},\ldots,t_{i_k})$ denote the Alexander
polynomial of $L_{\xi}$ (where $L_\xi=l_{i_1}\cup\cdots\cup
l_{i_k}$).

For the trivial representation  $\mathbf 1\!: G\to \mathbf C^*$,
$L_{\mathbf 1} =\emptyset$, and we set $\Delta_{L_{\mathbf 1}}=1$.
Let $\hat G^{(1)}$ be the subset of representations $\xi \in\hat G$
such that $L_{\xi}$ consists of a single component:
$L_\xi=l_{i(\xi)}$.

Finally,  $\vert H_1(\hat M^3_{\pi},\mathbf Z)\vert$ denotes the
cardinality of $H_1(\hat M^3_{\pi},\mathbf Z)$ when its finite, or
zero when it is infinite. The extension of Mayberry-Murasugi's
formula to homology spheres is the following.

\begin{theorem}
\label{thm:mainthm}
 In the situation described above we have:
\[
\vert H_1(\hat M^3_{\pi},\mathbf Z)\vert=\pm\prod_{\xi\in\hat G}
 \Delta_{L_{\xi}}(\xi(m_{i_1}),\ldots,\xi(m_{i_k}))
\ \frac {\vert G\vert}{\prod\limits_{\xi\in\hat G^{(1)}}
(1-\xi(m_{i(\xi)}))}.
\]
\end{theorem}

The relationship between the Alexander polynomial and Franz-Reidemeister
torsion was discovered by Milnor in \cite{MilnorAP} and further developed
by Turaev \cite{TuraevKT}, who provided new proofs for classical results.
In particular \cite{TuraevKT} reproved Fox's formula for knots in homology
spheres, but not Mayberry-Murasugi's, which was said to require additional
considerations going beyond the scope of the paper \cite{TuraevKT}.

\subsection*{Acknowledgements} I am indebted to M. Sakuma for an
encouraging conversation.

\section{Franz-Reidemeister Torsion}

We review the basic notions and results about Franz-Reidemeister torsion
needed in this paper. See \cite{MilnorWT} and \cite{TuraevKT} for details.

\subsection{Torsion of a chain complex.} Let $F$ be a field and
$C_n\overset\partial\to C_{n-1}\overset\partial\to\cdots\overset\partial
\to C_0$ a chain complex of finite dimensional $F$-vector spaces. Choose
$c_i$ a basis for $C_i$ and $h_i$ a basis for the $i$-th homology group.
We shall define the torsion of $C_i$ with respect to those basis.

Choose $b_i$ a basis for the $i$-dimensional boundary space (the
image of $\partial: C_{i+1}\to C_i$) and a lift $\tilde b_i$,
which is a subset of $C_{i+1}$ such that $\partial \tilde
b_i=b_i$. It is easy to check that the union $b_i\cup h_i\cup
\tilde b_{i-1}$ is a basis for $C_i$. Let $[b_ih_i\tilde
b_{i-1}/c_i]\in F^*$ denote the determinant of the transition
matrix between both basis (its entries are the coordinates of
vectors in $b_i\cup h_i\cup \tilde b_{i_1}$ with respect to
$c_i$). We define:
\[
\tau (C_*;c_i,h_i)=\prod_{i=0}^n [b_ih_i\tilde
b_{i-1}/c_i]^{(-1)^{i+1}}\in F^*/\{\pm 1\}.
\]
 It can be checked that
this torsion is independent of the choice of the $b_i$ and it is well
defined up to sign. In addtion, if we change the basis $c_i$ and $h_i$ we
get:
\begin{equation}
\label{eqn:changebasis}
 \tau (C_*;c_i',h_i')=\tau
(C_*;c_i,h_i)\prod_{i=0}^n\left(\frac{[h_i'/h_i]}{[c_i'/c_i]}\right)^{(-1)^{i+1}}.
\end{equation}
Notice that we follow the convention of \cite{MilnorAP} and
\cite{TuraevKT} for the sign $(-1)^{i+1}$, opposite to the one of
\cite{MilnorWT}.

\subsection{Torsion of a cell complex.} Let $K$ be a finite $CW$-complex
and $\varphi\! :\pi_1K\to F^*$ a representation. We define the
complex with coefficients twisted by $\varphi$
\[
C_*(K;\rho)=C_*(\tilde K;\mathbf Z)\otimes_{\varphi}F,
\]
where $C_*(\tilde K;\mathbf Z)$ is the complex with integer
coefficients on the universal covering. When $\varphi=\mathbf 1$
is the trivial representation, $C_*(K;\mathbf 1)=C_*(K;F)$ is the
usual untwisted complex.

We now choose a \emph{canonical basis} for $C_i(K;\varphi)$, that will
play the role of $c_i$ in the definition of torsion. Let
$\{e^i_1,\ldots,e^i_{j(i)}\}$ be the $i$-dimensional cells of $K$. We lift
them to the universal covering and we take $c_i=\{\tilde e^i_1\otimes
1,\ldots,\tilde e^i_{j(i)}\otimes 1\}$. The basis $c_i$ is called a
\emph{canonical basis}. Choosing again a basis $h_i$ for the homology we
define:
\[
\tau(K;\varphi,h_i)=\tau(C_*(K;\varphi);c_i,h_i)\in
F^*/\pm\varphi(\pi_1K).
\]
This definition only depends on the combinatorial class of $K$,
 $\varphi$ and the $h_i$.

\begin{remark}
We add the indeterminacy $\varphi(\pi_1K)$ due to the choice of the lift
of cells $\tilde e^i_j$. Turaev avoids this indeterminacy by using Euler
structures.
\end{remark}

\begin{example}
\label{ex:homolgy}
 Let $N^n$ be an $n$-dimensional rational homology sphere, so that
$H_i(N^n;\mathbf Z)$ is finite for $1\leq i\leq n-1$. Let
$h_{n}\in H_n(N^n,\mathbf Z)$ denote the fundamental class and
$h_0$ a generator for $H_0(N^n,\mathbf Z)$. For the trivial
representation $\mathbf 1:\pi_1N^n\to \mathbf C$ we have
\cite{TuraevKT}:
\[
\tau( K^n;\mathbf 1, h_0, h_n)=\pm \prod_{i=1}^{n-2} \vert
H_i(N^n;\mathbf Z) \vert ^{(-1)^{i+1}}.
\]
\end{example}

\begin{example} Let $L$ be a link in an 3-dimensional integer homology sphere $M^3$ with $\mu$
components. Consider its exterior $E(L)=M^3-N(L)$.
 We view the group ring $\mathbf C[\mathbf Z^{\mu}]$ as
the Laurent polynomial ring with $\mu$ variables $\mathbf
C[t_1^{\pm 1},\ldots,t_{\mu}^{\pm 1}]$ and let $F=\mathbf
C(t_1,\ldots,t_{\mu})$ denote its fraction field. Consider the
representation induced by abelianization $\rho\!: \pi_1E(L)\to
H_1(E(L);\mathbf Z)\cong\mathbf Z^{\mu}\hookrightarrow F$. Suppose
that the Alexander polynomial of the link $\Delta_L$ is non-zero
(which is always the case for a knot). Then $C_*(E(L);\rho)$ is
acyclic and:
\[
\tau(E(L);\rho)=\left\{
\begin{array}{cl}
    \Delta_L(t_1,\ldots,t_{\mu}) & \textrm{ when } \mu>1,\\
    \frac{\Delta(t_1)}{t_1-1} & \textrm{ when }\mu=1.
    \end{array}
\right.
\]
This was proved by Milnor when $\mu=1$ and Turaev when $\mu>1$,
\cite{MilnorAP,TuraevKT}. Notice that those identities hold true
up to multiplication by a factor $\pm t_1^{\alpha_1}\cdots
t_{\mu}^{\alpha_\mu}$. The complex $C_*(E(L);\rho)$ has
non-trivial first homology precisely when $\Delta_L=0$. So the
formula holds true if we define $\tau(E(L);\rho)=0$ when the
complex is not acyclic.
\end{example}

\subsection{The order of a module over a Noetherian UFD}
Both examples above can be deduced from a theorem of Turaev, as both rings
$\mathbf Z$ and $\mathbf C[t_1^{\pm 1},\ldots,t_{\mu}^{\pm 1}]$ are
Noetherian unique factorization domains (UFD).

Let $R$ be a Noetherian UFD and $D$ a finite generate $R$-module.
The module $D$ has a presentation matrix, with $m$ rows and $n$
columns, where $n$ is the rank. We can always assume $m\geq n$, by
adjoining rows of zeros if necessary. The \emph{elementary ideal}
is the ideal generated by the minors of the presentation matrix of
size $n$, and the \emph{order} of $D$ is the greatest common
divisor of this elementary ideal. We denote it by $\vert D\vert$.

For instance when $R=\mathbf Z$, $\vert D\vert$ is the cardinality
of $D$ when finite or $0$ when infinite. For a link in a homology
sphere, $\Delta_L=\vert H_1(\widetilde {E(L)},\mathbf Z)\vert$,
where $\widetilde {E(L)}$ is the maximal abelian covering of the
exterior of the link.

\begin{theorem}[\cite{TuraevKT}]
Let $C_*$ be  a complex of free $R$-modules and let $F$ be the
fraction field of $R$. Then $C_*\otimes_R F $ is acyclic iff
$\vert H_i(C)\vert \neq 0$, $\forall i=0,\ldots,n$. In this case:
\[
\tau(C_*;c_i)=\prod_{i=0}^m\vert H_i(C_*)\vert^{(-1)^{i+1}}.
\]
\end{theorem}

Notice that  for a link exterior $\vert
H_0(\widetilde{E(L)},\mathbf Z)\vert=(t-1)$ when $\mu=1$, and
\linebreak
 $\vert H_0(\widetilde{E(L)},\mathbf Z)\vert=1$
when $\mu>1$.

\section{Decomposing the $G$-complex $C_*(\hat K, \mathbf C)$}
\label{sec:decomposing}

Some of the material of sections~\ref{sec:decomposing} and
\ref{sec:nonac} is contained in \cite{Sakuma}. In particular most
of the results are contained there, but we give them again for
completeness and for fixing notation.

We use the notation of the introduction. Choose $K$ a CW-complex
such that $\vert K\vert =M^3$ and $L$ is a subset of the
$1$-skeleton. Let $\hat K$ be the induced CW-decomposition of
$\hat M_{\pi}^3$. By Example~\ref{ex:homolgy}, \( \vert H_1(\hat
M_{\pi};\mathbf Z)\vert = \tau ( \hat K;\mathbf 1, \hat h_0, \hat
h_3) \) where $\hat h_0$ and $\hat h_3$ are $ \mathbf Z$-basis for
$H_0(\hat M_{\pi}^3;\mathbf Z)$ and $H_3(\hat M_{\pi}^3;\mathbf
Z)$ respectively. Hence we want to study the chain complex $
C_*(\hat K;\mathbf 1)=C_*(\hat K;\mathbf C) $.

When the group ring $\mathbf C[G]$ is viewed as $G$-module, it
decomposes as a direct sum according to its representations:
\[
\mathbf C[G]=\bigoplus_{\xi\in\hat G}\mathbf C[{\mathbf f}_\xi],
\]
where ${\mathbf f}_\xi=\frac1{\vert G\vert}\sum\limits_{g\in G}
\xi(g^{-1}) g \in \mathbf C[G]$ (cf. \cite{Serre}). The element
${\mathbf f}_{\xi}\neq 0$ satisfies ${\mathbf f}_{\xi}^2={\mathbf
f}_{\xi}$ and $g\, {\mathbf f}_{\xi}=\xi(g)\, {\mathbf f}_{\xi}$.
Thus $\mathbf C[{\mathbf f}_\xi]$ is a one dimensional $\mathbf C
$-vector space, isomorphic to the $G$-module associated to $\xi\!:
G\to\mathbf C^*$.

The group $G$ acts naturally on the complex $C_*(\hat K;\mathbf C)$, thus
we have a decomposition of chain complexes:
\begin{equation}
\label{eqn:decomp}
 C_*(\hat K;\mathbf C)=
 \bigoplus_{\xi\in\hat G} {\mathbf f}_\xi \,C_*(\hat K;\mathbf C).
\end{equation}
Next we identify each subcomplex ${\mathbf f}_\xi \, C_*(\hat
K;\mathbf C)$, starting with the trivial representation $\mathbf
1$.

\begin{lemma}
\label{lem:isof1}
 There is a natural isomorphism ${\mathbf f}_{\mathbf 1}\, C_*(\hat K;\mathbf
C) \cong C_*(K;\mathbf 1) = C_*(K;\mathbf C)$.
\end{lemma}

\begin{proof}
We have a natural projection $C_*(\hat K;\mathbf C)  \to
C_*(K;\mathbf C)$ that restricts to  \linebreak ${{\mathbf
f}_{\mathbf 1}}\,C_*(\hat K;\mathbf C)  \to C_*(K;\mathbf C)$. To
construct its inverse, we map a chain $c\in C_*(K;\mathbf C)$ to
${\mathbf f}_{\mathbf 1} \hat c $, where $\hat c$ is any lift of
$c$. Since multiplication by ${\mathbf f}_{\mathbf 1}=\frac1{\vert
G\vert}\sum g$ is an average, this construction does not depend on
the lift and it is easily checked to be the inverse.
\end{proof}

Since the isomorphism of Lemma~\ref{lem:isof1} is natural, it induces an
isomorphsim in homology. Combining it with decomposition
(\ref{eqn:decomp}) we get:

\begin{corollary}
\label{cor:acyclicity}
 The covering $\hat M^3_{\pi}$ is a rational homology sphere iff $
{\mathbf f}_\xi C_*(\hat K;\mathbf C) $ has trivial first homology
group for every $\xi\in \hat G$, $ \xi\neq\mathbf 1$.
\end{corollary}

We view $K-L_{\xi}$ as a cell decomposition of the pair
$(E(L_\xi),\partial E(L_\xi))$ so that $\widetilde{K-L_{\xi}}$ is
a cell decomposition of $(\widetilde {E(L_{\xi})},\partial
\widetilde {E(L_{\xi})})$. The representation $\xi: G\to  \mathbf
C^*$ induces a representation $\pi_1 E(L_{\xi} )\to\mathbf C^*$,
also denoted by $\xi$, so that we can consider the complex:
\[
C_*(K-L_{\xi};\xi) =C_*(\widetilde{K-L_{\xi}};\mathbf
Z)\otimes_{\xi}\mathbf C.
\]

\begin{lemma}
\label{lem:isofxi} The complex $C_*(K-L_{\xi};\xi) $
 is naturally isomorphic to
${\mathbf f}_\xi \,C_*(\hat K;\mathbf C) $.
\end{lemma}

\begin{proof}
The  projection $\widetilde{K-L_{\xi}}\to \hat K$ induces a
natural map $C_*(K-L_{\xi};\xi)\to {\mathbf f}_\xi \,C_*(\hat
K;\mathbf C) $. It is straightforward to check that it is well
defined. Before constructing the inverse, notice that if $\hat
e^i_j$ is a cell of $\hat K$ that projects to $L_{\xi}$, then
there exists a $g\in G$ (the image of its meridian) such that $g\,
\hat e^{i}_j=\hat e^i_j$ and $\xi(g)\neq 1$. Thus ${\mathbf
f}_{\xi}\,\hat e^{i}_j ={\mathbf f}_{\xi} \,g\, \hat e^{i}_j=
\xi(g)\,{\mathbf f}_{\xi}\,\hat e^{i}_j=0$. This shows that we can
construct a map just by taking lifts of chains, which is easily
checked to be the inverse.
\end{proof}

\section{The non-acyclic case}
\label{sec:nonac}

\begin{lemma}
\label{lem:samehomology}
 The homology of $\mathbf{f}_{\xi} C_*(\hat K;\mathbf C)$ is
isomorphic to $H_*(E(L_\xi);\xi)$.
\end{lemma}

\begin{proof} By Lemma~\ref{lem:isofxi},
 the homology of the complex $\mathbf{f}_{\xi} C_*(\hat K;\mathbf C)$ is
 isomorphic to
 $H_*(E(L_\xi),\partial E(L_\xi);\xi)$. Using the exact sequence of the
 pair, it suffices to prove that $H_*(\partial E(L_\xi);\xi)=0$. Notice that
$\partial E(L_\xi)$ is a union of 2-dimensional tori, such that
the restriction of $\xi$ to each component is nontrivial. This
implies that $H^0(\partial E(L_\xi);\xi)=0$ because for each
component, the 0-cohomology group gives the subspace invariant by
the representation. A standard argument using duality and the
Euler characteristic proves the claim.
\end{proof}

\begin{proposition}
\label{prop:acyclicity}
 For $\xi\in\hat G$,
\[
H_1(E(L_\xi);\xi)=0\quad\textrm{ iff }\quad
\Delta_{L_\xi}(\xi(m_{i_1}),\ldots,\xi(m_{i_k}))\neq 0.
\]

\end{proposition}

\begin{proof}
Consider the evaluation map $\epsilon_\xi\!: \mathbf
C[t_{i_1}^{\pm 1}, \ldots, t_{i_k}^{\pm 1}]\to \mathbf C$, i.e.
\[
\epsilon_{\xi}(p(t_{i_1},\ldots,t_{i_k})) =
p(\xi(m_{i_1}),\ldots,\xi(m_{i_k})).
\]
 The short exact sequence
\[
0\to\ker\epsilon_\xi\to \mathbf C[t_{i_1}^{\pm 1}, \ldots,
t_{i_k}^{\pm 1}]\overset{\epsilon_\xi}\longrightarrow \mathbf C\to
0
\]
induces a long exact sequence in homology. A direct computation
shows that $H_0(E(L_{\xi});\mathbf C)=0$ and $
H_0(E(L_{\xi});\mathbf C[t_{i_1}^{\pm 1}, \ldots, t_{i_k}^{\pm
1}])\cong H_0(E(L_{\xi}); \ker\epsilon_\xi)\cong\mathbf C$. Thus
 we have a surjection
\begin{equation}
\label{eqn:surject}
H_1(E(L_\xi);\mathbf C[t_{i_1}^{\pm 1}, \ldots,
t_{i_k}^{\pm 1}])\to H_1(E(L_\xi);\xi)\to 0
\end{equation}
Since $\Delta_{L_\xi}$ is the order of $H_1(E(L_\xi);\mathbf
C[t_{i_1}^{\pm 1}, \ldots, t_{i_k}^{\pm 1}]) $, the map
(\ref{eqn:surject}) is zero iff
$\Delta_{L_\xi}(\xi(m_{i_1}),\ldots,\xi(m_{i_k}))\neq 0$.
\end{proof}

The following corollary is obtained in \cite{Sakuma}, where a
formula for the first Betti number is given. Here it follows from
Corollary~\ref{cor:acyclicity}, Lemma~\ref{lem:samehomology} and
Proposition~\ref{prop:acyclicity}, and proves the non-acyclic case
of Theorem~\ref{thm:mainthm}:

\begin{corollary}
The covering $\hat M^3_\pi$ is a homology sphere iff
$\tilde\Delta_L(\xi)\neq 0$ for all $\xi\in \hat G$.
\end{corollary}

\section{Proof of the main theorem}

The strategy of the proof is as follows. The order of $H_1(\hat
M^3_{\pi};\mathbf Z)$ is the torsion of the complex $C_*(\hat
K;\mathbf C)$. We use the decomposition (\ref{eqn:decomp}) to
write this torsion as product of torsions  of the complexes
$\mathbf f_{\xi} C_*(\hat K;\mathbf C)$ (Formula~\ref{eqn:product}
below). To get this formula, we change the canonical basis for
$C_*(\hat K;\mathbf C)$ to a union of canonical basis for $\mathbf
f_{\xi} C_*(\hat K;\mathbf C)$ (this is done in
Subsection~\ref{subsec:changing}). In
Subsection~\ref{subsec:tosioncomplexes} we compute the torsion of
each complex in terms of Alexander polynomials. All computations
in Subsection~\ref{subsec:tosioncomplexes} have an indeterminacy
of roots of unit, since the torsions we compute are defined up to
some root of unity and the Alexander polynomial is defined up to
some factor $t_1^{\alpha_1}\cdots t_\mu^{\alpha_\mu}$. This
indeterminacy  is discussed in
Subsection~\ref{subsection:indeterminacy}.

\subsection{Changing the canonical basis.}\label{subsec:changing}

Let $\{e^i_j\mid i=0,1,2,3\textrm{ and } j=1,\ldots,j(i)\}$ denote
the set of cells of $K$. Choose lifts $\hat e^i_j$ to $\hat K$, so
that, for $i=0,1,2,3$,
\[
c_i =\{g\, \hat e^i_j \mid  j=1,\ldots,j(i)\textrm{ and } g\in
G/Stab(\hat e^i_j)\}
\]
is the set of $i$-cells of $\hat K$, and hence a canonical basis
for $C_i(\hat K;\mathbf C)$. Define
\[
c_i'(\xi)=\{{\mathbf f}_{\xi}\,\hat e^i_j \mid j=1,\ldots,
j(i)\textrm{ and } \xi \vert_{Stab(\hat e^i_j)}\textrm{ is
trivial}\}.
\]

\begin{lemma}
\label{lem:naturalbasis}
 The isomorphism  of Lemma~\ref{lem:isof1} maps $c'_i(\mathbf 1)$ to a
 canonical basis for $C_i(K;\mathbf C)=C_i(K; \mathbf 1)$. The one of Lemma~\ref{lem:isofxi}
maps $c_i'(\xi)$ to a canonical basis for $C_i(K-L_{\xi};\xi)$.
\end{lemma}

\begin{proof}
A direct computation shows that ${\mathbf f}_{\xi}\,\hat e^i_j$ is mapped
to $\tilde e^i_j\otimes_{\xi} 1$ when $\xi\neq\mathbf 1$, and  to $e^i_j$
when $\xi=\mathbf 1$. Counting elements, we realize that this is a
canonical basis.
\end{proof}

In particular $c'_i(\xi)$ is a basis for $\mathbf f_{\xi} C_i(\hat
K;\mathbf C)$ and  $\bigcup\limits_{\xi\in\hat G}c_i'(\xi)$ is a
basis for $C_i(\hat K; \mathbf C)$.

\begin{lemma}
\label{lem:candc'} $\prod\limits_{i=0}^3[\bigcup\limits_{\xi\in\hat
G}c_i'(\xi)/c_i]^{(-1)^{i}}=1$.
\end{lemma}

\begin{proof} For each subgroup $H< G$, the set of lifts $\{\hat
e^i_j\}$ of cells that have precisely $H$ as stabilizer, has zero Euler
characteristic. This implies that there are cancellations in the
alternated product.
\end{proof}

When computing the torsion of the complexes $C_i(K;\mathbf C)$ and
$C_i(K-L_{\xi};\xi)$, we will assume that we are using the
canonical basis of Lemma~\ref{lem:naturalbasis}.

 It follows from decomposition (\ref{eqn:decomp}) and from Lemmas~\ref{lem:isof1}, \ref{lem:isofxi} and \ref{lem:candc'}
  that
\begin{equation}
\label{eqn:product}
 \tau (\hat K;\mathbf 1, \hat h_0, \hat h_3)=
 \tau( M^3;\mathbf 1, \hat h_0,\hat h_3)
  \prod _{\begin{array}{c}\scriptstyle{\xi\in\hat
G}\\ \scriptstyle{\xi\neq\mathbf 1}\end{array}}\tau
(E(L_{\xi}),\partial E(L_{\xi});\xi) .
\end{equation}

\subsection{The torsion $\tau (E(L_{\xi}),\partial E(L_{\xi});\xi)$ as evaluation of the Alexander polynomial.}
\label{subsec:tosioncomplexes}

\begin{lemma}  $\tau(M^3;\mathbf 1 ,\hat h_3,\hat h_0 )=\vert
G\vert$.
\end{lemma}

\begin{proof} Since the isomorphism of Lemma~\ref{lem:isof1}  is induced
by the projection $\hat M_\pi^3\to M^3$, $\hat h_0$ is mapped to
$h_0$, a generator for $H_0(M^3;\mathbb Z)$, and the fundamental
class $\hat h_3$ of $\hat M^3_\pi$ is mapped to $\vert G\vert$
times the fundamental class $h_3$ of $M^3$. Thus, by
(\ref{eqn:changebasis}):
\[
\tau( M^3;\mathbf 1,\hat h_3,\hat h_0)=\vert G\vert\, \tau(
M^3;\mathbf 1,   h_3,  h_0)
\]
and $\tau( M^3;\mathbf 1,  h_3,  h_0)=1$ because $M^3$ is a
homology sphere (Ex.~\ref{ex:homolgy}).
\end{proof}

\begin{lemma}
 \label{lem:tauxi}
 Let $\xi\in\hat G$ with $\xi\neq\mathbf 1$. Then
\[
\tau(E(L_{\xi}),\partial E(L_{\xi});\xi)=\left\{
\begin{array}{ll}
\tilde \Delta_{L}(\xi)& \textrm{if }\xi\not\in \hat G^{(1)},\\
 \frac{\tilde \Delta_{L}(\xi)}{\xi(m_{i(\xi)})-1} &  \textrm{if }\xi \in \hat
 G^{(1)}.
\end{array}
\right.
\]
\end{lemma}

\begin{proof} We consider the exact sequence of the pair
\[
0\to C_*(\partial E(L_\xi);\xi)\to C_*(E(L_\xi);\xi) \to
C_*(E(L_\xi),\partial E(L_\xi);\xi)\to 0
\]
We showed in the proof of Lemma~\ref{lem:samehomology} that
$C_*(\partial E(L_\xi);\xi)$ is acyclic. By
Sections~\ref{sec:decomposing} and \ref{sec:nonac}, we may assume
that all complexes in the sequence are acyclic. Applying
\cite[Thm.~3.1]{MilnorWT} (see also \cite[Lemma~4]{MilnorCD}), we
get
\[
\tau(E(L_{\xi});\xi)=\tau(\partial E(L_{\xi}),\xi) \  \tau
(E(L_{\xi}),\partial E(L_{\xi});\xi).
\]
The torsion $\tau(\partial E(L_{\xi}),\xi) $ is trivial because $\partial
E(L_{\xi})$ is even dimensional \cite{Franz} (see also \cite{MilnorAP},
otherwise a direct computation on the torus shows it). Thus
\[
 \tau (E(L_{\xi}),\partial E(L_{\xi});\xi) = \tau(E(L_{\xi});\xi).
\]

Consider the representation $\rho\! :\pi_1E(L_{\xi})\to \mathbf
C[t_{j_1}^{\pm 1},\ldots t_{j_k}^{\pm 1}] \subset \mathbf
C(t_{j_1},\ldots t_{j_k})$ corresponding to the abelianization, so
that $H_*(E(L_{\xi});\rho)=0$ and
\[
\tau(E(L_{\xi});\rho)= \left\{
\begin{array}{cl}
    \Delta_L(t_{j_1},\ldots,t_{j_k}) & \textrm{ when }
    \xi\not\in\hat G^{(1)} \textrm{ (i.e. } k>1 \textrm{)}, \\
    \frac{\Delta(t_{i(\xi)})}{t_{i(\xi)}-1} & \textrm{ when } \xi\in\hat
    G^{(1)}.
    \end{array}
\right.
\]
Representations $\xi$ and $\rho$ are related by the evaluation morphism:
\[
\begin{array}{rcl}
    \epsilon_{\xi}\! : \mathbf C[t_{j_1}^{\pm 1},\ldots t_{j_k}^{\pm
1}] & \to & \mathbf C\\
f(t_{j_1},\ldots, t_{j_k}) & \mapsto & p(\xi(m_{j_1}),\ldots,
\xi(m_{j_k})).
\end{array}
\]
We have $\xi=\epsilon_{\xi}\circ\rho$. We claim  that
\begin{equation}
\label{eqn:naturalness}
 \epsilon_{\xi}(\tau(E(L_\xi);\rho))=\tau
(E(L_{\xi});\xi).
\end{equation}
Before proving the claim, it is relevant to notice that
$\epsilon_\xi$ is defined on the polynomial ring $\mathbf
C[t_{j_1}^{\pm 1},\ldots t_{j_k}^{\pm 1}]$ but not on the whole
fraction field $\mathbf C(t_{j_1},\ldots, t_{j_k})$. This problem
can be avoided by computing the torsion following the method of
 \cite{TuraevKT} or \cite{MilnorAP}. Namely, we choose $\tilde b_2$ to be the
canonical basis for $C_3(E(L_{\xi});\xi)$. Acyclicity implies that
$b_2=\partial \tilde b_2$ is a set of linearly independent
elements in $C_2(E(L_{\xi});\xi)$. We complete $b_2$ to a basis
for $C_2(E(L_{\xi});\xi)$ by choosing elements of the canonical
basis, whose union we denote by $\tilde b_1$. Again acyclicity
implies that $b_1=\partial \tilde b_1$ is a set of linearly
independent elements, and so on. Each time we choose elements of
the canonical basis for $\xi$, we do the corresponding choice for
$\rho$. In this way, all the determinants involved in the torsion
for $\rho$ belong to $\mathbf C[t_{j_1}^{\pm 1},\ldots
t_{j_k}^{\pm 1}]$ and have the property that $\epsilon_\xi$ maps
them to the corresponding determinants for computing the torsion
for $\xi$. Hence (\ref{eqn:naturalness}) follows.
\end{proof}

\subsection{Eliminating the indeterminacy of roots of unit}
\label{subsection:indeterminacy}

The argument above proves the Theorem~\ref{thm:mainthm} up to some
indeterminacy corresponding to roots of unit (appearing in
Lemma~\ref{lem:tauxi}).
We claim that the formula holds true without
this indeterminacy if we use the \emph{same choice of the Alexander
polynomial for each sublink.} Given a representation $\xi\in\hat G$,
its complex conjugate $\overline\xi$ is also a representation in
$\hat G$. In addition $L_{\xi}=L_{\overline\xi}$, hence
$\Delta_{L_{\xi}}=\Delta_{L_{\overline\xi}}$ and
$\Delta_{L_{\overline\xi}}
(\overline\xi(m_{i_1}),\ldots,\overline\xi(m_{i_k}))$ is the complex
conjugate of $ \Delta_{L_{\xi}}(\xi(m_{i_1}),\ldots,\xi(m_{i_k}))$.
Thus the right hand term of the formula in Theorem~\ref{thm:mainthm}
belongs to $\mathbf R$ and the claim is proved.

\section{Generalizations}

\subsection{Abelian coverings branched along graphs}

Following \cite{Sakuma}, we can generalize
Theorem~\ref{thm:mainthm}
 to coverings of homology three-spheres along graphs. There are
some important restrictions to our graph. Firstly, the vertices
must have valency three, if we want that the covering is a
manifold. In addition, the only abelian finite subgroups of
$SO(3)$ are either cyclic or $\mathbf Z/2\mathbf Z\oplus \mathbf
Z/2\mathbf Z$. Hence, when we have a trivalent vertex, the
ramification on the adjacent edges must be $2$. Once those
restrictions are established, Theorem~\ref{thm:mainthm}
generalizes for $L$ to be such an embedded graph, not only a link.
This is because each cyclic subgroup of $SO(3)$ fixes an edge, and
even if $L$ is a graph, for every representation $\xi\in \hat G$,
the subgraph $L_{\xi}$ is a link.

\subsection{Higher dimensional knots}

We work in the PL-category. Let $M^{n+2}$ a $(n+2)$-dimensional
homology sphere and $K^n\subset M^{n+2}$ an $n$-knot, with
$K^n\cong S^n$. We can consider the $d$-cyclic branched covering
$\hat M^{n+2}_d\to M^{n+2}$ branched along $K^n$.

In this case, we do not have only a the Alexander polynomial but
several Alexander invariants. The exterior
$E(K^n)=M^{n+2}-N(K^{n})$ has the homology of the circle and we
consider its infinite abelian covering $\widetilde{E(K^n)}$.

The $i$-th Alexander invariant is defined to be the order  of
$H_i(\widetilde{E(K^n)}; \mathbf Z)$ as $ \mathbf
Z[t,t^{-1}]$-module:
\[
A_i(t)=\vert H_i(\widetilde{E(K^n)}; \mathbf Z)\vert.
\]

\begin{theorem} The covering $M^{n+2}_d$ is a rational homology
sphere iff $A_i(\zeta)\neq 0$ for every $d$-root of unit $\zeta$
and every $i$. When it is a rational homology sphere
\[
\prod_{i=1}^n \vert H_i(\hat M^{n+2}_d; \mathbf Z)\vert^{(-1)^{i+1}}=
\prod_{i=1}^{n+1}\prod_{k=1}^d A_i(\zeta^d)^{(-1)^{i+1}},
\]
where $\zeta$ is a primitive d-root of unit.
\end{theorem}

Of course this formula is only relevant when $n$ is odd, because when $n$
is even each one of the products is $1$ (the torsion of an even
dimensional manifold is trivial \cite{Franz}).

The proof of this theorem follows exactly the same argument as
Theorem~\ref{thm:mainthm} with minor changes.

\end{document}